\documentclass[11pt]{article}
\usepackage{authblk} 
\usepackage[margin=1in]{geometry}
\usepackage{times}  % DO NOT CHANGE THIS
\usepackage{helvet}  % DO NOT CHANGE THIS
\usepackage{courier}  % DO NOT CHANGE THIS
\usepackage[hyphens]{url}  % DO NOT CHANGE THIS
\usepackage{graphicx} % DO NOT CHANGE THIS
\usepackage{natbib}  % DO NOT CHANGE THIS AND DO NOT ADD ANY OPTIONS TO IT
\usepackage{caption} % DO NOT CHANGE THIS AND DO NOT ADD ANY OPTIONS TO IT
\usepackage{amsthm,amssymb,amsmath}
\newtheorem{definition}{Definition}
\newtheorem{lemma}{Lemma}
\newtheorem{theorem}{Theorem}
\newtheorem{axiom}{Axiom}
\usepackage{coqdoc}
\usepackage[hidelinks]{hyperref}
\title{A Topological Rewriting of Tarski's Mereogeometry}

\author[1]{Patrick Barlatier}
\author[1]{Richard Dapoigny}

\affil[1]{University of Savoie Mont-Blanc, LISTIC}
\affil[ ]{\texttt{patrick.barlatier@univ-smb.fr}}
\affil[ ]{\texttt{richard.dapoigny@univ-smb.fr}}

\date{}
\begin{document}
	\maketitle
	\begin{abstract}Qualitative spatial models based on Goodman-style mereology and pseudo-topology often pose problems for advanced geometric reasoning, as they lack true Euclidean geometry and fully developed topological spaces. We address this issue by extending an existing formalization, grounded in an underlying type theory using the \textit{Coq} language, together with the Whitehead-like point-free Tarski's geometry. More precisely, we leverage an available library called $\lambda$-MM to formalize Tarski’s geometry of solids by investigating an algebraic formulation of topological relations on top of the library. Given that Tarski’s work is grounded in Le{\'s}niewski’s mereology, and despite the fact that $\lambda$-MM barely implements Tarski's geometry, the first part of the paper supplements their work by proving that mereological classes correspond to regular open sets. It forms a topology of individual names extensible with Tarski’s geometric primitives. Unlike classical approaches used in qualitative logical theories, we adopt a solution that enables the specification of a topological space from mereology and a geometric subspace, thereby enhancing the theory’s expressiveness. Then, in a second part, we prove that Tarski’s geometry forms a subspace of the previous topology in which regions are restricted classes. We prove three postulates of Tarski’s work reducing his axiomatic system and extend the theory with the T2 (Hausdorff) property and additional definitions.
	\end{abstract}
\begin{center}
\begin{minipage}{0.85\textwidth}
\footnotesize
\textbf{Note.} This is the full version of the paper accepted at AAAI-26.
The arXiv version includes the complete list of authors.
\end{minipage}
\end{center}
\vspace{1em}

\begin{center}
\begin{minipage}{0.85\textwidth}
\url{https://www.univ-smb.fr/listic/technologies/logiciels/modeles-decidables-pour-la-specification-d-ontologies/extending-tarskis-mereo-geometry-with-topological-definitions/}
\end{minipage}
\end{center}

\section{Introduction}\label{intro}
	In the field of Qualitative Spatial Reasoning (QSR), the interplay between mereology and topology is a common foundation for many theories \cite{Biacino1991,Roeper1997,Duntsch2004,Hahmann2013}. Theories resulting from the interplay between mereology and topology follow the ideas of Whitehead’s point-free geometry \cite{Whitehead1929} in which points are substituted by regions of space. Building on mereology either using a quasi-topological theory or by introducing topological notions to increase its expressiveness, e.g., in contact algebras or in RCC8, has gained significant attention in QSR research \cite{Menger1940,Varzi1994,Smith1996,Mormann2001}. However, these approaches often face challenges such as (i) the formal specification of topological boundaries (i.e., what kind of entity a boundary is, either a spatial object or an abstract entity?), (ii) the decidability of the resulting systems (the question of determining tractable fragments is still unanswered for many theories) \cite{Hahmann2012,Hahmann2013}, (iii) insufficient reasoning mechanisms without which the spatial representations are less useful and (iv) a shallow representation of mereology using a direct mapping of part of relations in first-order logic and resulting in poor expressiveness. To address these problems, we rely on a previous work, i.e., the $\lambda$-MM open source library \cite{Barlat2023} which specifies an algebraic structure grounded in a syntactic, type-theoretic representation of mereology implemented in the Coq theorem prover. The central objective is (i) to extend this library with an algebraic formulation of topological relations using appropriate primitives as an alternative to contact-based mereotopologies and (ii) to provide a unified theory incorporating mereology, geometry and topology. For these purpose, we will extend the $\lambda$-MM library implemented as a syntactic model grounded in \textit{Coq} \citep{Bertot2004,Coquand1988}. The library includes both mereology and Tarski’s geometry of solids within an inductive type space, all while maintaining the classical logic required by many mereological lemmas. In this paper, we first show in section \ref{foundations} by adding some definitions, that the boolean mereological model implemented in $\lambda$-MM can be seen as a structure of regular open sets. In section \ref{TM}, we extend the mereological structure with topological operators (interior, boundary, closure) and demonstrate that it satisfies Kuratowski’s axioms. Finally, in section \ref{four} we introduce the notion of saturated interior point, and after redefining topological operators on a sub-base of the previous topology, we prove Tarski’s axioms and the $T2$ property.
	\section{Foundations of $\lambda$-MM}\label{foundations}
	The $\lambda$-MM library heavily relies on Le{\'s}niewski's stratified systems known as Le{\'s}niewski's Ontology (LO) and Mereology. 
	\subsection{From Le{\'s}niewski's Ontology to $\lambda$-MM}
	The Logic of Names (LO) is a formal system that manipulates names as symbolic terms referring to individuals, collections, or aggregates. It includes rules of symbolic manipulation comparable to those of formal logic. Like predicate calculus, LO is a general framework for logical inference, but it uses a grammar centered on names instead of predicates and quantifiers. It is structured around three basic syntactic categories: names (terms), propositions (statements), and functors (constructors that derive new names or propositions from existing ones). At the core of LO lies the functor $\eta$, which acts as a binary syntactic relation between names and produces well-formed sentences. The term $\eta \: A \: x$ asserts that $A$ is a name falling under the classification designated by the name $x$. This relation is not set-theoretic membership, but a purely internal syntactic relation defined by the structural composition of names. For example, $\exists \: B, \eta \: B \: A$ states that if $B$ refers to the neighbor’s daughter”, hence $B$ designates the same entity as “Nicole” ($A$), without relying on an external model. Unlike classical logic which relies on a metalanguage and interpretation over an external domain, LO formulates such assertions within the language itself, based on syntactic inference. The formal role of $\eta$ is similar to the ”:” symbol in type theory \cite{Simon1987}, where $x \: : \: A$ asserts that a term belongs to a given type. However, type theory assumes a strict stratification between objects and types, often distributed across multiple layers of meta-theory. By contrast, the sentence $\eta \: x \: a$ in LO operates within a single unified syntactic layer, where names, categories, and relations are treated homogeneously. Above these assumptions, $\eta$ replaces external typing hierarchies with internal derivation rules, removing the need for separate semantic interpretation.
	\par In $\lambda$-MM, the dependent type theory underlying Coq provides a constructive and syntactic model of LO. Names are objects of an inductive type $N$, defined by a grammar that encodes their internal structure. The type system specifies not only how names are formed, but also how their composition and well-formedness are guaranteed by construction. This replaces the classical extension/intension distinction with a purely syntactic logic where (i) the existence of a name is guaranteed by its inductive construction, (ii) identity is determined syntactically, not semantically and (iii) the internal structure of a name determines its relations.
	\par The functor $\eta$ is interpreted in $\lambda$-MM as a judgmental relation of type $N \rightarrow N \rightarrow Prop$. A sentence $\eta \: A \: b$ means that $A$ satisfies the inductive specification of $b$, and this judgment is derivable from the internal grammar of $b$. As such, $\eta$ no longer relies on extensional semantics, but on proof-theoretic rules internal to the system.
	\par In $\lambda$-MM, we distinguish between simple relations, of type $N \rightarrow N \rightarrow Prop$, and compound relations, defined as compositions of $\eta$ with unary or binary relations $\phi$. A compound relation $R = \eta \: \circ \: \phi$  has type $N \rightarrow Prop$, where $\phi : N \rightarrow N$ or $\phi : (N \rightarrow N) \rightarrow N$. The following lemma holds for such relations:
	\begin{lemma}\label{In_to_eta} 
		$\forall \: a : N , \: \forall \: X : object , \;  X \: \in \:(\phi \: a) \: \leftrightarrow \: \eta \: (\iota X) \: (\phi \: a)$	
	\end{lemma}
	in which $\iota$ denote the singleton constructor\footnote{It maps every object $X$ to its associated name (set), $\iota \: X$.}. A list of relations is given in table \ref{table1}. Compound relations are translated into a set of simple relations. For example, the name conjunction is defined as: $$\forall a \: b, a \cap b \: \triangleq \: \forall P, \; (\eta \: P \: a \wedge \eta \: P \: b)$$
	\begin{table}[!h]
		\begin{tabular}[!h]{ l }
			\begin{tabular}[!h]{ c }
				\textbf{Simple relations} \\	
				\begin{tabular}{ c | l | l |}	
					$\eta$ & singular inclusion & $\forall P \; Q, \eta \: P \: Q \;$ \\
					$\equiv$ & equality (individual names) & $\forall P \; Q, P \: \equiv \: Q $ \\
					$\approx$ & equality (plural names) & $\forall a \: b, a \approx b$ \\
					$\subseteq$ & name inclusion & $\forall a \: b, a \subseteq b$ \\			
				\end{tabular} \\
			\end{tabular} \\
			\begin{tabular}[!h]{ c }		
				\textbf{Compound relations} \\
				\begin{tabular}{ c | l | l |}
					$\mathit{neg}$ & negation of name & $\forall P \; Q, \eta \: P \; (neg \; Q)$ \\
					$\cap$ & name conjunction & $\forall a \: b, a \cap b$  \\
					$\cup$ & name disjunction & $\forall a \: b, a \cup b$  \\
				\end{tabular}
			\end{tabular}
		\end{tabular}
		\caption[Table 1]{Basic relations in LO between names in $\lambda$-MM.}\label{table1}
	\end{table}
	Two constant plural names are specified in LO, which are respectively the empty name ($\Lambda$) defined as the contradictory name, and the universal name ($V$).
	\subsection{From Le{\'s}niewski's Mereology to $\lambda$-MM}
	Le{\'s}niewski's mereology redefines the notion of a class not as a set of members, but as a concrete whole (a mereological sum) composed of parts. These wholes, referred to as mereological classes or m-classes, are treated as individual entities within the system. In mereology, (i) an m-class or aggregate is itself a single entity and (ii) an element is understood as a part of the whole, not a member in the set-theoretic sense. The primitive relation is no longer membership, but the part–whole relation in the following sense: writing $\eta \: A \: klass \: a$ means that $A$ now refers to a totality (all parts recursively inside $a$), whereas writing $\eta \: A \: a$ rather assumes that $A$ only refers to a single member of the "$a$'s (like membership). All entities (both m-classes and their parts) share the same formal status: they are all names of individuals within a nominalist system.
	\par It is well established that m--classes are extensional, meaning they are fully determined by their parts. As highlighted in \cite{Asenjo1977}, this type of extensional class contrasts with intensional or structured aggregates, which depend on how their constituents are organized. Hence, m-classes fall under the category of extensional classes, as opposed to intensional classes where the mode of composition affects identity.
	\par In the mereological component of $\lambda$-MM, all relations are compound: using $\phi$ as the generic name defined above, all mereological constructors such as $\mathit{pt}$ (part-of), $\mathit{ppt}$ (proper-part-of), $\mathit{klass}$ (m--class), etc. satisfy the semantics of lemma \ref{In_to_eta}. These functions operate within the inductive type $N$ introduced earlier. As an example, the definition of a m--class formalizes the composition of $\mathit{klass}$ with $\eta$:
	\begin{definition}$\mathtt{(m-class)}$\label{mklass}\mbox{}\\
		$\forall A \: a, \: \eta \: A \: (klass \: a) \triangleq (Individual \: A \: \wedge$\\ $(\forall \: B, \eta \: B \: a \rightarrow \eta \: B \: (pt \: A)) \: \wedge$\\
		$(\forall \: B, \eta \: B (pt \: A) \rightarrow \exists \: C \: D, \eta \: C \: a \wedge \eta \: D (pt \: C) \wedge \eta \: D (pt \: B)))$.
	\end{definition}
	The term $Individual \: A$ is equivalent to the term $\eta \: A \: A$ which guarantees that $A$ is an individual name.
	\par In \cite{Barlat2023}, the authors formally demonstrate that a full implementation of mereology, comprising: (i) a lower layer based on the Calculus of Inductive Constructions (as implemented in Coq, restricted to classical logic via appropriate libraries), (ii) the $LO^1$ fragment of the Logic of Names, limited to monadic second-order primitives and (iii) a core mereological sub-theory centered on the part-of primitive ($\mathit{pt}$) yields a model that is isomorphic to a quasi-Boolean algebra\footnote{A Boolean algebra without a bottom element.} (see e.g., \cite{Tarski1956b,Clay1974,Barlat2023}) under the assumption that the partial order coincides with the $\mathit{pt}$ relation.
	\subsection{Extending $\lambda$-MM with Topological Notions} 
	In Tarski's geometry of solids \cite{Tarski1956a} which is based on mereology, a solid is defined as a sum of spheres a.k.a. balls. This sum is not simply an arbitrary aggregation but is defined by internal geometric properties, derived from a primitive relation (the ball) and geometric relations such as tangency, concentricity, mereological inclusion, etc. Thus, the identity of a solid depends on how it is constructed from balls, and not solely on the extensional set of balls it contains. This dependence on a construction rule (and not on extension alone) is proof of an intensional class.
	\par Solids that are put into one-to-one correspondence with classes of interior points form regular open sets (axioms 2 and 3 in \cite{Tarski1956a}). This relation of inwardness itself relies on a mereological structure (Def. 9), which encodes a dependence between balls and the points they contain. Consequently, a solid is characterized not only by its interior points, but by the balls from which it is generated, yielding a clear status of geometric intension. Tarski explicitly states (Def. 6) that certain concepts, such as that of a point, are defined by a class of concentric balls. Such a definition highlights the existence of a hierarchy of constructions in which entities are generated by relations, not just by their extension. Hence, points are already intensional classes (classes of balls), and bodies, as sums of geometrically determined balls, are even more clearly intensional mereological classes. 
	\par The formalization in Coq of a pointless geometry, inspired by Tarski and based on Le{\'s}niewski's mereology, makes it possible to build a library with strong theoretical and practical potential, particularly for modeling qualitative, topological, and constructive spaces. The objectives of the present work are (i) to replace primitive points with mereological classes of balls, (ii) to construct topological regions as regular open sets (e.g., from sums of balls), (iii) to specify spatial relations (inclusion, neighborhood, boundary) in mereological terms, thus without resorting to membership or sets and (iv) to provide a formal alternative to spatial ontologies like RCC8 or set-theoretic geometries.
	\section{A Topological Interpretation of Mereology}\label{TM}
	\subsection{Introduction of meet and join operators}
	In $\lambda$-MM, the whole space of mereology is defined as an individual name called $Universe$ such as:
	\begin{definition}$\mathtt{(Universe)}$\label{Univ}\mbox{}\\
		$\forall \: P \: : N, \; \eta \: P \: Universe \triangleq \eta \: P \: (klass \: V).$
	\end{definition}
	The Boolean model underlying $\lambda$-MM is structured using the part-of relation ($\mathit{pt}$) as its ordering principle, but lacks meet and join operators. Since their role matters for the topological interpretation of mereology, we first add their definitions to $\lambda$-MM. The Boolean model admits a least upper bound (the supremum, or join operator). In this setting, an individual name can never be empty, whereas a plural name may be. Since no bottom element exists, a general meet operator cannot be defined globally, though restricted meet operations can still be introduced locally whenever defined. We have extended the initial version of $\lambda$-MM while preserving the model’s partial structure. The Boolean join of two objects $P$ and $Q$ is represented by an object $R$, defined as the nominal disjunction of their respective parts, in accordance with the underlying nominalist mereology:
	\begin{definition}$\mathtt{(join)}$\label{join}\mbox{}\\
		$\forall P \: Q \: R : N, \; \eta \: R \: (b\_sum \: P \: Q) \; \triangleq \; (\eta \: P \: P \; \wedge \; \eta \: Q \: Q \; \wedge$\\
		$  \eta \: R \: (klass \: ((pt \: P) \cup (pt \: Q)))).$
	\end{definition}
	Similarly, the meet $R$ simply results from the nominal conjunction of their parts.
	\begin{definition}$\mathtt{(meet)}$\label{meet}\mbox{}\\
		$\forall P \: Q \: R : N, \; \eta \: R \: (b\_prod \: P \: Q) \; \triangleq \; (\eta \: P \: P \; \wedge \; \eta \: Q \: Q \; \wedge$\\
		$ \eta \: R \: (klass \: ((pt \: P) \cap (pt \: Q)))).$
	\end{definition}
	Then, it is routine to prove that these operators satisfy the standard properties such as (i) uniqueness, (ii) an individual and its complement are disjoint and (iii) their sum is the universe. 
	\subsection{Mereology as a Topology of Individual Names}
	Mereological classes play the role of non-empty open subsets of some topological structure $\langle N, \mathbb{T}_{MM}\rangle$. The basis of the topology $\mathbb{T}_{MM}$ is a collection of names from $N$ restricted to singular names (i.e., $Ind$). For open sets the interior of $Q$ is precisely $Q$, which is asserted both with a first term asserting that $Q$ is an individual and a second term which returns a value identical to $Q$. Hence, the interior is an individual, i.e., a m--class:
	\begin{definition}$\mathtt{(interior)}$\label{interior}\mbox{}\\
		$\forall \: P \: Q \: : N, \; \eta \: P \: (interior \: Q) \triangleq \: \eta \: Q \: Q \: \wedge \: \eta \: P \: Q$
	\end{definition}
	By adding the $interior$ definition to the extension of $\lambda$-MM we can prove that the resulting structure is a topology $\langle N, \mathbb{T}_{MM}\rangle$, over a basis of individual names. We assume following \cite{Mormann2001}, that the set of mereological individuals in $\lambda$-MM is interpreted as the set of open subsets of $\mathbb{T}_{MM}$ in a one-to-one correspondence. In such a way, topological relations like inclusion, intersection, union and complement are respectively and faithfully related to mereological operations such as $pt$, $b\_product$, $b\_sum$ and $compl$. Here, $compl$ is the mereological complement such that $P$ is the complement of $Q$ relative to $R$ ($\eta \: P \: (relCompl \: Q  \: R)$) in which $R$ is substituted by $Universe$). 
	Instead of relying on the classical Kuratowski closure axioms, we adopt their interior dual formulation \cite{Pervin1964}. This choice simplifies the proof strategy in our setting, where interiors correspond to individual names, and open sets are mereological classes.
	\begin{enumerate}
		\item Preservation of the total space: $int(\mathbb{T}_{MM}) = \mathbb{T}_{MM}$.
		\item Intensiveness: $\forall A \subseteq \mathbb{T}_{MM}: int(A) \subseteq A$.
		\item Idempotence: $\forall A \subseteq \mathbb{T}_{MM}: int(int A)) = (int  A)$.
		\item Preservation of binary intersections: $\forall A, B \subseteq \mathbb{T}_{MM}: int(A \cap B) = int(A) \cap int(B)$
	\end{enumerate}
	\par Notice that if subsets of $\mathbb{T}_{MM}$ are all open, the second axiom reduces to $int(A) = A$. Following Tarski’s second postulate which identifies the set of all interior points of a region as a regular open set, we formalize this idea using a topological class that encodes Kuratowski’s interior axioms within our nominal and inductive framework. Mathematical structures can be defined in Coq as type classes with fields defining basic data structures and other fields describing the properties that the basic data must satisfy. Alternatively, basic data can be transmitted as arguments to the type class (parameterized type class). A major interest of type classes is that they are generic in the sense where they are parameterized with data, and changing data will not affect the structure of the type class. The topology is described as a type class in which each field defines the type of fundamental properties corresponding to the interior operator over the family of subsets of $\mathbb{T}$ as described in definition \ref{TopoTypeClass}. Each field corresponding to Kuratowski interior axioms depends here on the fact that individual names are open sets (the $Open$ predicate assumes that $Q$ is an individual name). The fields labeled $open\_XXX$ stand for interior theorems and all depend on the $Interior$ specification:
	\begin{definition}$\mathtt{(Topological \: class)}$\label{TopoTypeClass}\mbox{}\\
		\begin{coqdoccode}
			\coqdockw{Class} \coqdocvar{Kuratowski\_Open} (\coqdocvar{X} :\coqdockw{Type})(\coqdocvar{Univ} :\coqdocvar{X})\coqdoceol
			\coqdocindent{3.00em}(\coqdocvar{equiv\_indiv} : \coqdocvar{X} \ensuremath{\rightarrow} \coqdocvar{X} \ensuremath{\rightarrow} \coqdockw{Prop})\coqdoceol
			\coqdocindent{3.00em}(\coqdocvar{Prod} : \coqdocvar{X} \ensuremath{\rightarrow} \coqdocvar{X} \ensuremath{\rightarrow} \coqdocvar{X}) := \{\coqdoceol
			\coqdocindent{1.00em}
			\coqdocvar{Open} \coqdocindent{5.00em} : \coqdocvar{X} \ensuremath{\rightarrow} \coqdockw{Prop} ;\coqdoceol
			\coqdocindent{1.00em}
			\coqdocvar{Interior} \coqdocindent{4.00em}: \coqdocvar{X} \ensuremath{\rightarrow} \coqdocvar{X};\coqdoceol
			\coqdocindent{1.00em}
			\coqdocvar{open\_space} \coqdocindent{2.30em}: \coqdocvar{equiv\_indiv} (\coqdocvar{Interior} \coqdocvar{Univ}) \coqdocvar{Univ}; \coqdoceol
			\coqdocindent{1.00em}
			\coqdocvar{open\_intensive} \coqdocindent{1.00em}: \ensuremath{\forall} \coqdocvar{Q} :\coqdocvar{X}, \coqdocvar{Open} \coqdocvar{Q} \ensuremath{\rightarrow} \coqdoceol
			\coqdocindent{8.00em}(\coqdocvar{equiv\_indiv} (\coqdocvar{Interior} \coqdocvar{Q}) \coqdocvar{Q}); \coqdoceol
			\coqdocindent{1.00em}
			\coqdocvar{open\_idempotent}             : \coqdockw{\ensuremath{\forall}} \coqdocvar{Q} :\coqdocvar{X}, \coqdocvar{Open} \coqdocvar{Q} \ensuremath{\rightarrow} (\coqdocvar{equiv\_indiv} \coqdoceol
			\coqdocindent{8.00em}(\coqdocvar{Interior} (\coqdocvar{Interior} \coqdocvar{Q})) (\coqdocvar{Interior} \coqdocvar{Q}));\coqdoceol
			\coqdocindent{1.00em}
			\coqdocvar{open\_inter} \coqdocindent{2.80em} : \coqdockw{\ensuremath{\forall}} \coqdocvar{A} \coqdocvar{B} : \coqdocvar{X}, \coqdocvar{Open} \coqdocvar{A} \ensuremath{\land} \coqdocvar{Open} \coqdocvar{B} \ensuremath{\rightarrow}\coqdoceol
			\coqdocindent{8.00em} \ensuremath{\forall} \coqdocvar{P} \coqdocvar{Q} \coqdocvar{R}:\coqdocvar{X},  (\coqdocvar{equiv\_indiv} \coqdocvar{P} (\coqdocvar{Interior} \coqdocvar{A})) \coqdoceol
			\coqdocindent{8.00em}\ensuremath{\land} (\coqdocvar{equiv\_indiv} \coqdocvar{Q} (\coqdocvar{Interior} \coqdocvar{B})) \ensuremath{\land} \coqdoceol
			\coqdocindent{8.00em}
			(\coqdocvar{equiv\_indiv} \coqdocvar{R} (\coqdocvar{Prod} \coqdocvar{A} \coqdocvar{B})) \ensuremath{\rightarrow}\coqdoceol
			\coqdocindent{8.00em} \coqdocvar{equiv\_indiv} (\coqdocvar{Prod} \coqdocvar{P} \coqdocvar{Q}) (\coqdocvar{Interior} \coqdocvar{R}); \coqdoceol
			\coqdocindent{1.00em}\}.	
		\end{coqdoccode}	
	\end{definition}
	Instances are defined in the Coq system through the \textit{Instance} keyword. For this instance to be created, we need to prove that these fields are populated with appropriate theorems and we have to prove each of them. It is obtained using a Coq instance of the type class by providing the corresponding theorems (for each field, the instance match the corresponding theorem name):
	\begin{theorem}
		\begin{coqdoccode}
			\coqdocemptyline
			\coqdocnoindent\coqdockw{Instance} \coqdocvar{proof\_topo} : \coqdocvar{Kuratowski\_Open} \coqdocvar{N}\coqdoceol
			\coqdocindent{9.00em} \coqdocvar{Universe} \coqdocvar{singular\_eq} \coqdocvar{b\_product} := \{\coqdoceol
			\coqdocindent{1.00em}
			\coqdocvar{Open} \coqdocvar{Q} \coqdocindent{5.00em} := $\eta$ \coqdocvar{Q} \coqdocvar{Q};\coqdoceol
			\coqdocindent{1.00em}
			\coqdocvar{Interior} \coqdocvar{Q} \coqdocindent{4.10em} := \coqdocvar{interior} \coqdocvar{Q}; \coqdoceol
			\coqdocindent{1.00em}
			\coqdocvar{open\_space} \coqdocindent{3.40em} := \coqdocvar{int\_space};\coqdoceol
			\coqdocindent{1.00em}
			\coqdocvar{open\_intensive} \coqdocindent{2.10em} := \coqdocvar{int\_open};\coqdoceol
			\coqdocindent{1.00em}
			\coqdocvar{open\_idempotent} \coqdocindent{1.2em} := \coqdocvar{int\_idempotent};\coqdoceol
			\coqdocindent{1.00em}
			\coqdocvar{open\_inter} \coqdocindent{3.80em} := \coqdocvar{int\_to\_intersections};\coqdoceol
			\coqdocnoindent
			\}.\coqdoceol	
		\end{coqdoccode}
	\end{theorem}
	Corresponding proofs are given in the link as \textit{Coq} theorems. They have been proved in a semi-automated way using CoqHammer \cite{Czajka2018} (mostly with \textit{Vampire} and \textit{Eprover}). Since open sets refer to individuals and provided that the interior is open, the closure is introduced as a standard topological construction, i.e., with the complement of the interior of the complement. Formally, it is introduced as follows:
	\begin{definition}\label{}$\mathtt{(closure)}$\label{closure}\mbox{}\\
		$\forall \: P \: Q, \: \eta \: P \: (closure \: Q) \triangleq \; (\eta \: Q \: Q \wedge \; (exists \: R \: S,$\\ $\eta \: R \: (compl \: Q) \: \wedge \: \eta \: S \: (interior \: R) \: \wedge \: \eta \: P \: (compl \: S))).$	
	\end{definition}
	The first term assumes that $R$, as an individual, is identified with its unique complement and characterizes an open set. The second term follows the same rule, since the interior is unique, and the resulting individual $S$ refers to an open set. Then, regularity is derived with asserting that the interior of the closure of $Q$ coincides with $Q$:
	\begin{lemma}$\mathtt{(regular)}$\label{regul}\mbox{}
		\begin{coqdoccode}
			\ensuremath{\forall} \coqdocvar{P} \coqdocvar{Q}, (\coqdocdefinition{$\eta$} \coqdocvar{Q} \coqdocvar{Q} \ensuremath{\land} \coqdocdefinition{$\eta$} \coqdocvar{P} (\coqdocvar{closure} \coqdocvar{Q}) \ensuremath{\rightarrow} \coqdocvar{Q} $\equiv$ (\coqdocvar{interior} \coqdocvar{P})).\coqdoceol
		\end{coqdoccode}	
	\end{lemma}
	The domain of mereology, i.e., the class of individuals, is interpreted as the class of regular open sets without the empty set (it does not mean that the resulting lattice is atomless). The mereological boundary of a name $Q$ is then defined as the intersection of the closure of $Q$ and the closure of its complement:
	\begin{definition}$\mathtt{(boundaryM)}$\label{bound0}\mbox{}\\
		$\forall \: P \: Q \: : N, \: \eta \: P \: (boundaryM \: Q) \: \triangleq \: (\eta \: Q \: Q \; \wedge$ \\ $\eta \: P \: ((closure \: Q) \: \cap \: (closure \: (compl \: Q))).$
	\end{definition}			
	Using definition \ref{bound0}, we prove that this boundary is empty establishing that the topology is clopen:
	\begin{theorem}$\mathtt{(clopen)}$\label{clopen}\mbox{}
		\begin{coqdoccode}
			\ensuremath{\forall} \coqdocvar{Q}, $\eta$ \coqdocvar{Q} \coqdocvar{Q} \ensuremath{\land} \ensuremath{\lnot} (\coqdocvar{Q} $\equiv$ \coqdocvar{Universe}) \ensuremath{\rightarrow} \coqdocvar{boundaryM} \coqdocvar{Q} $\approx$ $\Lambda$.\coqdoceol
		\end{coqdoccode}	
	\end{theorem}
	While mereology can be regarded, with the addition of a few definitions, as a topology based on individual names, its limitation to representing clopen sets is too weak to allow for effective applications. By contrast, its specification within Tarski’s geometric context offers significant potential for spatial reasoning.
	\section{From Tarski's geometry to Topology}\label{four}
	We now connect the formal system $\lambda$-MM with Tarski’s geometry of solids\footnote{We rather consider that "solid" or "spatial regions" will be replaced with general term "region".} using balls as primitives. This allows us to reconstruct Tarski’s postulates using purely topological and mereological primitives. Mereology and regular open sets play a fundamental role in Tarski’s Foundations of Geometry of Solids \cite{Tarski1929}. His system has been shown to be semantically complete with respect to Euclidean models of $\mathbb{R}^n$ \cite{Bennett2001}. Following \cite{Vakarelov2007}, we adopt the interpretation of regular open sets of a topological space as models for mereotopology. Since Tarski’s system integrates mereological reasoning, we can reconstruct Euclidean geometry by merging geometrical and mereological features within a topological extension of $\lambda$-MM. We adopt the view of \cite{Borgo2010} in which points are defined as filters of concentric balls rather than primitives. This allows a purely constructive and nominalist reformulation of Tarski’s postulates using Coq to extend our topological extension of $\mathbb{T}_{MM}$ with geometrical definitions inside the basis of a more complete topological version. In the next section, we use mereology supplemented with topological definitions for the specification of geometrical relations.
	\subsection{Topological Interpretation of Tarski's Geometry of Solids}\label{TopoTarski}
	Tarski has given a model for atomless Boolean algebra, consisting of the family of regular open sets of an Euclidean space and the relation of set-inclusion \cite{Tarski1956b}. Partly based on mereology, Tarski has suggested a categorical axiomatization of a pointless\footnote{It does not mean that we have no points in the theory but rather that points do not appear among the primitive notions and by contrast are defined by means of other primitives.} geometry of regions using only a ball (or sphere) as primitive \cite{Tarski1929}. This seminal work has given rise to many applications in region-based theories such as \cite{Bennett2001,Borgo2010,Gruszcz2008,Dapo2015} among others. The link with "the method of extensive abstraction" developed by A. Whitehead focuses on the specification of equidistance between points (points are seen as infinite sets of concentric balls). Tarski used the fact that these points are neither elements of, nor external to a given ball (this idea has also been taken up by Ja{\'s}kowski \cite{Jaskowski1948}). He implies thereby that a boundary exists between two individuals. Once mereology is proved to be a space for regular open sets, we suggest here to express adequate definitions for the geometry of regions close to Tarski’s geometry of solids, but using topological primitives such as interior and boundary. The fact that singular names correspond regular open sets, is crucial to understand postulates $2$ and $3$ of Tarski's paper \cite{Tarski1929} which formalize the connection between the structure (solid, pt) and (regular open set, inclusion). Furthermore, we have the ability to redefine Tarski’s definitions with topological primitives based on points as sets of concentric balls. We specify a topology in which the sub-base consists in the sum of balls whereas open sets are individuals. Assuming that we have \textit{balls} as a primitive name included into $N$, we may start out with a subspace of $V$ which is characterized as the m--class of all open balls:
	\begin{definition}$\mathtt{(Gspace)}$\label{gspace}\mbox{}\\
		$\forall \: P, \: \eta \: P \: Gspace \: \triangleq \: \eta \: P \: (klass \: balls).$	
	\end{definition}
	Since $Gspace$ is an individual, the m--class of all balls corresponds to an open set. We assume that each ball refers to an open ball. The following axiom states that there exists at least a ball in $Gspace$:
	\begin{axiom}$\mathtt{(ball\_exists)}$\label{allBalls}\mbox{}
		\begin{coqdoccode}
			$! \: balls$.	
		\end{coqdoccode}	
	\end{axiom}
	We assume also a one-to-one correspondence between $Gspace$ and the usual 3D \textit{Euclidian} space leading to an interpretation of primitives in the class of regular regions of $\mathbb{R}^3$. From these assumptions, we are able to prove the existence and uniqueness of $Gspace$, and that any ball is a part of $Gspace$.
	\par Introducing a set of definitions, such as external ($ET$) and internal tangency ($IT$), externally diametral ($ED$) and internally diametral tangencies ($IT$) (see e.g., \cite{Tarski1929} for details), we reach the major definition of concentricity ($Concent$). It results from the above-mentioned set of definitions. It shows up all properties of an equivalence relation (reflexive, symmetric and transitive). Finally a point is defined as the set of all balls all of them being concentric with a given one.
	\begin{definition}$\mathtt{(Point)}$\label{point}\mbox{}\\
		$\forall \: P \: Q : N, \: \eta \: P \: (Point  \: Q) \; \triangleq \; (\eta \: P \: balls \: \wedge \: \eta \: Q \: balls \: \wedge$ \\ $\eta \: P \: (Concent \: Q)).$ 
	\end{definition}
	Notice that $(\mathit{Concent} \: Q)$ is a plural, that is, a set of all balls concentric with ball $Q$. While balls denote first-order objects, points as set of balls are second-order objects and fall within the scope of a monadic second-order theory. From the point definition, two lemmas about points are worth detailing, i.e., concentric points and reflexivity.
	\begin{lemma}$\mathtt{(equiv\_points)}$\label{eqPoints}\mbox{}
		\begin{coqdoccode}
			\coqdockw{\ensuremath{\forall}} \coqdocvar{P} \coqdocvar{Q} \coqdocvar{R}, \coqdocdefinition{$\eta$} \coqdocvar{P} (\coqdocvar{Point} \coqdocvar{R}) \ensuremath{\land}\coqdoceol \coqdocindent{2.00em}  \coqdocdefinition{$\eta$} \coqdocvar{Q} (\coqdocvar{Point} \coqdocvar{R}) \ensuremath{\rightarrow}  \coqdocdefinition{$\eta$} \coqdocvar{P} (\coqdocvar{Concent} \coqdocvar{Q}).		
		\end{coqdoccode}	
	\end{lemma}
	\begin{lemma}$\mathtt{(Point\_refl)}$\label{ptRefl}\mbox{}
		\begin{coqdoccode}
			\coqdockw{\ensuremath{\forall}} \coqdocvar{P}, \coqdocdefinition{$\eta$} \coqdocvar{P} \coqdocvar{balls} \ensuremath{\rightarrow} \coqdocdefinition{$\eta$} \coqdocvar{P} (\coqdocvar{Point} \coqdocvar{P}).		
		\end{coqdoccode}	
	\end{lemma}
	Another important geometric relations can be introduced, the betweenness between points. Given two points, a third point lies between them iff there exists three balls concentric with them which are in the external diametric relation:
	\begin{definition}$\mathtt{(Betweenness)}$\label{betw}\mbox{}\\
		$\forall A B C, \eta A \: (btw B \: C) \triangleq (\eta A \: balls \wedge \eta B \: balls \wedge \eta \: C \: balls \wedge$\\ $\exists \: D \: E \: F, \: \eta \: D \: (Point \: B) \wedge \eta \: E \: (Point \: A) \wedge \eta \: F \: (Point \: C) \wedge$\\ $\eta \: D (ED \: E \: F)).$
	\end{definition}
	Regions follow Tarski's definition as set of balls:
	\begin{definition}$\mathtt{(Region)}$\label{region}\mbox{}\\
		$\forall \: P \: : N, \: \eta \: P \: Region \triangleq (Individual \: P \wedge \exists \: b, \: b \subseteq balls \: \wedge$ \\ $\eta \: P \: (klass \: b)).$
	\end{definition}
	A similar definition has been given in \cite{Clay2021}. It is more precise than Tarski's description since it involves a particular subset of balls giving rise to a corresponding m--class. However, Tarski has pointed out that the concept of "sum" should be understood as collection in Le{\'s}niewski's formalism. Precisely, using definition \ref{region} we can derive the following equivalence which fully agrees with Tarski's view:
	\begin{lemma}$\mathtt{(region\_as\_coll)}$\label{regColl}\mbox{}
		\begin{coqdoccode}
			\coqdockw{\ensuremath{\forall}} \coqdocvar{P}, \coqdocdefinition{$\eta$} \coqdocvar{P} \coqdocvar{Region} \ensuremath{\leftrightarrow} \coqdoceol \coqdocindent{2.00em} \coqdocdefinition{$\eta$} \coqdocvar{P} (\coqdocvar{coll} \coqdocvar{balls}).	
		\end{coqdoccode}	
	\end{lemma}
	Since a region is a set of balls and provided that balls are open sets, their union is open as well. Lemmas \ref{ballreg}, \ref{spacereg} and \ref{reghasball} show respectively that any ball is a solid (the converse is not true), that the geometric space is also a solid and that any region includes balls.
	\begin{lemma}$\mathtt{(any\_ball\_is\_a\_region)}$\label{ballreg}\mbox{}
		\begin{coqdoccode}
			\coqdockw{\ensuremath{\forall}} \coqdocvar{P}, \coqdocdefinition{$\eta$} \coqdocvar{P} \coqdocvar{balls} \ensuremath{\rightarrow} \coqdoceol \coqdocindent{2.00em}\coqdocdefinition{$\eta$} \coqdocvar{P} \coqdocvar{Region}.
		\end{coqdoccode}
	\end{lemma}
	\begin{lemma}$\mathtt{(whole\_space\_is\_a\_region)}$\label{spacereg}\mbox{} 
		\begin{coqdoccode}
			\coqdoceol
			\coqdocindent{2.00em}\coqdockw{\ensuremath{\forall}} \coqdocvar{P}, \coqdocdefinition{$\eta$} \coqdocvar{P} \coqdocvar{Gspace} \ensuremath{\rightarrow} \coqdocdefinition{$\eta$} \coqdocvar{P} \coqdocvar{Region}.
		\end{coqdoccode}
	\end{lemma}
	\begin{lemma}$\mathtt{(any\_region\_includes\_balls)}$\label{reghasball}\mbox{}
		\begin{coqdoccode}
			\coqdoceol
			\coqdocindent{2.00em}	
			\coqdockw{\ensuremath{\forall}} \coqdocvar{P}, \coqdocdefinition{$\eta$} \coqdocvar{P} \coqdocvar{Region} \ensuremath{\rightarrow} \ensuremath{\exists} \coqdocvar{C}, \coqdocdefinition{$\eta$} \coqdocvar{C} \coqdocvar{balls} \ensuremath{\land} \coqdocdefinition{$\eta$} \coqdocvar{C} (\coqdocvar{pt} \coqdocvar{P}).
		\end{coqdoccode}
	\end{lemma}
	Skipping details, we now focus on the last definition proposed by Tarski. The point $a$ is an interior point of the solid $Q$ if there exists a ball $P'$ which is at the same time an element of the point $a$ and a part of the region $Q$:
	\begin{definition}$\mathtt{(InteriorPoint)}$\label{intpt}\mbox{}\\
		$\forall \: P \: Q, \: \eta \: P \: (InteriorPoint \: Q) \: \triangleq$ \\ $(\eta \: Q \: Region \: \wedge \: \exists \: P', \: \eta \: P \: (Point \: P') \: \wedge \: \eta \: P' \: (pt \: Q)).$
	\end{definition}
	When Tarski refers to the point $a$ he involves a set of all concentric balls to a given ball, therefore $Point \: P'$ plays here the role of $a$ and $P'$ is one of its balls, provided that ball $P'$ is a part of region $Q$. Thus, points are seen as filters of regions. In summary the definition of the interior point describes a set which in Le{\'s}niewski's formalism can be related to a m--class. It is important to see that such a definition asserts that while all the referred points belong to the region, their related concentric balls expand outside the region. This aspect is crucial for the definition of the topological interior of a region.
	\par Beside lemmas proving that for any region, there exists at least an interior point and that any individual point is a ball we can formalize Tarski's postulates 2 and 3 as theorems. Since an interior point corresponds to a set, the m--class of this set is an individual according to mereological principles hence, it refers to a regular open set:
	\begin{theorem}$\mathtt{(Tarski\_P2)}$\label{TP2}\mbox{}\\
		\begin{coqdoccode}
			\coqdockw{\ensuremath{\forall}} \coqdocvar{Q}, \coqdocdefinition{$\eta$} \coqdocvar{Q} \coqdocvar{Region} \ensuremath{\rightarrow} \ensuremath{\exists} \coqdocvar{R}, \coqdocdefinition{$\eta$} \coqdocvar{R} (\coqdocvar{klass} (\coqdocvar{InteriorPoint} \coqdocvar{Q})).
		\end{coqdoccode}
	\end{theorem}
	\begin{theorem}$\mathtt{(Tarski\_P3)}$\label{TP3}\mbox{} \\
		\begin{coqdoccode}
			\coqdockw{\ensuremath{\forall}} \coqdocvar{R} \coqdocvar{Q}, \coqdocdefinition{$\eta$} \coqdocvar{R} (\coqdocvar{klass} (\coqdocvar{InteriorPoint} \coqdocvar{Q})) \ensuremath{\rightarrow} \coqdocdefinition{$\eta$} \coqdocvar{R} \coqdocvar{Region}.		
		\end{coqdoccode}
	\end{theorem}
	Theorem \ref{TP2} states that if $Q$ is a region, the m--class of all its interior points is non-empty regular open set. This results from the properties of any m--class. Alternatively, in theorem \ref{TP3} if we have a class of points that is a regular open set (i.e., a m--class), then it corresponds to an individual $R$ which is a region. The interpretation of postulate 4 demonstrates that the inclusion of sets  of interior points for two regions implies that their respective open sets are in the part-of relation:
	\begin{theorem}$\mathtt{(Tarski\_P4)}$\label{TP4}\mbox{}
		\begin{coqdoccode}
			\coqdoceol
			\coqdocindent{2.00em}\ensuremath{\forall} \coqdocvar{P} \coqdocvar{Q} \coqdocvar{R} \coqdocvar{S}, \coqdocvar{$\eta$} \coqdocvar{P} \coqdocvar{Region} \ensuremath{\land} \coqdocdefinition{$\eta$} \coqdocvar{Q} \coqdocvar{Region} \ensuremath{\land}\coqdoceol
			\coqdocindent{2.00em} (\coqdocvar{InteriorPoint} \coqdocvar{P}) $\subseteq$ (\coqdocvar{InteriorPoint} \coqdocvar{Q}) \ensuremath{\land} \coqdoceol
			\coqdocindent{2.00em}\coqdocdefinition{$\eta$} \coqdocvar{R} (\coqdocvar{klass} (\coqdocvar{InteriorPoint} \coqdocvar{P})) \ensuremath{\land}\coqdoceol
			\coqdocindent{2.00em}
			\coqdocdefinition{$\eta$} \coqdocvar{S} (\coqdocvar{klass} (\coqdocvar{InteriorPoint} \coqdocvar{Q})) \ensuremath{\rightarrow} \coqdocdefinition{$\eta$} \coqdocvar{R} (\coqdocdefinition{pt} \coqdocvar{S}).
		\end{coqdoccode}
	\end{theorem}	
	\subsection{The Topology of Regions} As underlined in subsection \ref{TopoTarski}, the original specification of interior points is not sufficient to construct the base for a topology. Therefore, we have slightly modified Tarski's definition \ref{intpt} to constrain balls of interior points to be saturated, i.e., constrained to be inside the region:
	\begin{definition}$\mathtt{(SatInteriorPoint)}$\label{satintpt}\mbox{}\\
		$\forall \: P \: Q, \: \eta \: P \: (sat\_InteriorPoint \: Q) \triangleq (\eta \: Q \: Region \: \wedge$\\ $\exists \: P', \: \eta \: P \: (Point \: P') \: \wedge \: \eta \: P' \: (pt \: Q) \: \wedge \: (\eta \: P' \: (pt \: P) \: \rightarrow$\\ $(P \: \equiv \: P'))).$	
	\end{definition}
	All lemmas derived with Tarski's definition \ref{intpt} remain valid but we can prove important lemmas. For example, (i) the geometric space is an open set built with the m--class of its interior points and (ii) if $Q$ is an open region, then it coincides with the m--class of its interior points. Therefore, a simple formulation can emerge for the topological interior:
	\begin{definition}$\mathtt{(InteriorG)}$\label{intgeom}\mbox{}\\
		$\forall \: P \: Q \: : N, \: \eta \: P \: (InteriorG \: Q) \triangleq$ \\ $\eta \: P \: (klass \: (sat\_InteriorPoint \: Q)).$	
	\end{definition}
	that is, $P$ is the topological interior of any region $Q$ iff it is the m--class of its saturated interior points. Using the generic type class of definition (\ref{TopoTypeClass}), we are able to instantiate it with appropriate theorems to satisfy all Kuratowski's axioms:
	\begin{theorem}
		\begin{coqdoccode}
			\coqdocemptyline
			\coqdocnoindent \coqdockw{Instance} \coqdocvar{proof\_topo} : \coqdocvar{Kuratowski\_Open} \coqdocvar{N} \coqdocvar{Gspace} \coqdocvar{singular\_eq} \coqdocvar{b\_product} := \{\coqdoceol
			\coqdocindent{1.00em}
			\coqdocvar{Open} \coqdocvar{Q} \coqdocindent{5.00em} := $\eta$ \coqdocvar{Q} \coqdocvar{Region};\coqdoceol
			\coqdocindent{1.00em}
			\coqdocvar{Interior} \coqdocvar{Q} \coqdocindent{4.10em} := \coqdocvar{InteriorG} \coqdocvar{Q}; \coqdoceol
			\coqdocindent{1.00em}
			\coqdocvar{open\_space} \coqdocindent{3.20em}  := \coqdocvar{int\_Gspace};\coqdoceol
			\coqdocindent{1.00em}
			\coqdocvar{open\_intensive}   \coqdocindent{1.90em}  := \coqdocvar{int\_is\_open};\coqdoceol
			\coqdocindent{1.00em}
			\coqdocvar{open\_idempotent} \coqdocindent{0.80em}  := \coqdocvar{int\_pt\_idempotent};\coqdoceol
			\coqdocindent{1.00em}
			\coqdocvar{open\_inter}  \coqdocindent{3.50em} := \coqdocvar{int\_to\_intersections};\coqdoceol
			\coqdocnoindent	\}.\coqdoceol	
		\end{coqdoccode}
	\end{theorem}
	It follows that the geometric space with balls can be interpreted as a topology $\langle balls, \mathbb{T}_{MG}\rangle$, based on the \begin{coqdoccode}\coqdocvar{InteriorG}\end{coqdoccode} definition. To complete $\mathbb{T}_{MG}$, we first add the boundary definition already suggested by \cite{Jaskowski1948} and implicit in Tarski's definition of equidistance \cite{Tarski1929}. Informally, it says that a point $P$ belongs to the boundary of a given region $R$ iff $P$ is both not a part of $R$ and not external to $R$:
	\begin{definition}$\mathtt{(boundary)}$\label{boundary}\mbox{}\\
		$\forall \: P \: Q, \eta \: P \: (boundary \: Q) \triangleq (\eta \: Q \: Region \wedge \eta \: P \: balls \: \wedge$ \\ $\forall \: Y, \eta \: Y \: balls \: \wedge \: \eta \: P \: (Point \: Y) \rightarrow \neg \eta \: Y \: (pt \: Q) \wedge \neg \eta \: Y (ext \: Q)).$	
	\end{definition}
	The boundary definition assumes the existence of a plural name $boundary \: Q$ and ball $P$ is a member of this plural. Ball $P$ generates a point which is centered on the boundary. Using the boundary (denoted $\partial$), we can prove among others, that (i) interior and boundary points are disjoints, (ii) boundary is not a region or (iii) the boundary of the complement of a region coincides with the boundary of the region $\partial (complG \; P) \approx \partial \: P$.
	\begin{lemma}$\mathtt{(int\_not\_bound)}$\label{bounDisj}\mbox{}\\
		\begin{coqdoccode}
			\ensuremath{\forall} \coqdocvar{P} \coqdocvar{Q}, $\eta$ \coqdocvar{P} (\coqdocvar{sat\_InteriorPoint} \coqdocvar{Q}) \ensuremath{\rightarrow} \ensuremath{\lnot} $\eta$ \coqdocvar{P} ($\partial$ \coqdocvar{Q}).	
		\end{coqdoccode}
	\end{lemma}
	\begin{lemma}$\mathtt{(bound\_not\_reg)}$\label{not_reg}\mbox{}\\
		\begin{coqdoccode}
			\ensuremath{\forall} \coqdocvar{P} \coqdocvar{Q}, $\eta$ \coqdocvar{Q} \coqdocvar{Region} \ensuremath{\land} $\eta$ \coqdocvar{P} ($\partial$ \coqdocvar{Q}) \ensuremath{\rightarrow} $\neg$(($\partial$ \coqdocvar{Q}) $\approx$ \coqdocvar{Region}).\coqdoceol
		\end{coqdoccode}
	\end{lemma}	
	\begin{lemma}$\mathtt{(compl\_boundary)}$\label{compl_bound}\mbox{}\\
		\begin{coqdoccode}
			\ensuremath{\forall} \coqdocvar{P} \coqdocvar{Q}, $\eta$ \coqdocvar{P} (\coqdocvar{complG} \coqdocvar{Q}) \ensuremath{\rightarrow} $\partial$ \coqdocvar{Q} $\approx$ $\partial$ \coqdocvar{P}.
		\end{coqdoccode}
	\end{lemma}	
	By allowing both open and closed regions, a boundary of a region is also a boundary of its complement. In such a way we avoid the arbitrary choice of to which region the boundary actually belongs. In summary, the boundary is a set of points (concentric spheres) that do not belong mereologically either to a region nor to its complement. It bears some similarity with the discrete version of the RCC \cite{Roy2002} in that boundaries are defined as regions without interiors.
	\par The mereological complement $complG$ satisfies definition $relcompl$ (see section \ref{TM}) in which $R$ is substituted with $Gspace$. The topological complement is defined as the union of the mereological complement with the boundary. 
	\begin{definition}$\mathtt{(Tcompl)}$\label{Tcompl}\mbox{}\\
		$\forall \: P \: Q, \: \eta \: P \: (Tcompl \: Q) \: \triangleq \: \eta \: P \: ((complG \: Q) \: \cup \: (\partial \: Q)).$
	\end{definition}
	From the topological complement, the closure of $Q$ is the complement of the interior of the complement of $Q$:
	\begin{definition}$\mathtt{(closureG)}$\label{closG}\mbox{}\\
		$\forall \: P \: Q, \: \eta \: P \:  (closureG \; Q) \triangleq (\eta \: P \: Region \: \wedge \: \eta \: Q \: Region \: \wedge$ \\ $ \eta \: P \: (Tcompl \: (InteriorG \: (complG \: Q)))).$	
	\end{definition}
	From the closure definition we infer that the closure is inferentially equivalent to the union of the interior of $Q$ with its boundary as expected:
	\begin{theorem}$\mathtt{(clos\_to\_int\_bound)}$\label{closIntBound}\mbox{}\\
		\begin{coqdoccode}
			\ensuremath{\forall} \coqdocvar{P} \coqdocvar{Q},
			$\eta$ \coqdocvar{Q} \coqdocvar{Region} \ensuremath{\land} $\neg$(\coqdocvar{Q} $\equiv$ \coqdocvar{Gspace})
			\ensuremath{\rightarrow} \coqdoceol \coqdocindent{1.00em}$\eta$ \coqdocvar{P} (\coqdocvar{closureG} \coqdocvar{Q}) \ensuremath{\leftrightarrow}
			$\eta$ \coqdocvar{P} ((\coqdocvar{InteriorG} \coqdocvar{Q}) $\cup$ ($\partial$ \coqdocvar{Q}))
		\end{coqdoccode}
	\end{theorem}
	In Theorem \ref{closIntBound}, the definition of $closureG$ depends on complements, hence the complement of $Q$ must be non-empty ($Q \neq Gspace$). Using definition \ref{betw} (betweenness) and the basic geometric axiom stating that for any two distinct points there exists a third point between them, we can show that $Gspace$ satisfies the $\mathit{Hausdorff}$ property.
	\begin{theorem}$\mathtt{(is\_Hausdorff)}$\label{Hausdorff}\mbox{}\\
		\begin{coqdoccode}
			\ensuremath{\forall} \coqdocvar{P} \coqdocvar{Q}, $\neg$($\eta$ \coqdocvar{P} (\coqdocvar{Concent} \coqdocvar{Q})) \ensuremath{\rightarrow} \ensuremath{\exists} \coqdocvar{B1} \coqdocvar{B2}, $\eta$ \coqdocvar{B1} \coqdocvar{balls} \ensuremath{\land} \coqdoceol $\eta$ \coqdocvar{B2} \coqdocvar{balls} \ensuremath{\land} $\eta$ \coqdocvar{P} (\coqdocvar{sat\_InteriorPoint} \coqdocvar{B1}) \ensuremath{\land} \coqdoceol $\eta$ \coqdocvar{Q} (\coqdocvar{sat\_InteriorPoint} \coqdocvar{B2}) \ensuremath{\rightarrow} \ensuremath{\lnot}\ensuremath{\exists} \coqdocvar{Z}, \coqdoceol $\eta$ \coqdocvar{Z} (\coqdocvar{sat\_InteriorPoint} \coqdocvar{B1}) \ensuremath{\land} $\eta$ \coqdocvar{Z} (\coqdocvar{sat\_InteriorPoint} \coqdocvar{B2}).\coqdoceol
		\end{coqdoccode}
	\end{theorem}
	It follows that $Gspace$ is a $T2$ topology in which distinct points have disjoint neighborhoods. Furthermore, it complies with the assumption that $Gspace$ is isomorph to $\mathbb{R}^3$. Many definitions can be added to extend the expressiveness, such as the convexity relation:
	\begin{definition}$\mathtt{(convexRegion)}$\label{convex}\mbox{}\\
		$\forall \: P,  \: \eta  \: P  \: convexRegion  \; \triangleq  \; (\eta  \: P  \: Region  \: \wedge  \: \forall  \: Q  \: R \: S,$\\ $\eta \: Q \: (sat\_InteriorPoint \: P) \wedge \eta \: R (sat\_InteriorPoint \: P) \wedge$\\ $\eta \: S \: (btw \: Q \: R) \: \rightarrow  \: \eta \: S \: (sat\_InteriorPoint \: P)).$	
	\end{definition}
	Definition \ref{convex} intuitively asserts that given two interior points $Q$ and $R$ of a region $P$, then any point between $Q$ and $R$ must be an interior point of this region. It excludes situations in which the region has the shape of a "U" or a snake. 
	\par Let $\mathcal{M} = \langle R, \leq, \mathbb{T}_{GM} \rangle$ be the structure of Tarski’s geometry, where: $R$ is the set of individuals that are regions, i.e., a sub-base of $N$, $\leq$ is the mereological part-of relation and $\mathbb{T}_{GM}$ is the topology of regular open sets according to definition (\ref{TopoTypeClass}). The structure can be summarized as:
	\begin{enumerate}
		\item $ \mathcal{M} \models \mathit{Topological \; axioms}$, where \textit{Topological axioms} denotes Kuratowski's axioms for $\mathbb{T}_{GM}$.	
		\item $\mathcal{M} \models \mathit{Hausdorff}$, formally expressed as:
		\[\forall x,y \in balls \,(x \neq y \rightarrow \exists U,V \in \mathbb{T}_{GM}(x \in U \wedge y \in V \wedge\] 
		\[U \cap V = \emptyset)).\]	
		\item $\mathcal{M} \models \mathit{Convexity}$, which asserts that the open sets (bodies) satisfy a convexity condition:
		\[\forall R \in \mathbb{T}_{GM}, \forall x,y \in R, [x,y] \subseteq R.\]
		\item We may conjecture that the following homeomorphism holds:
		\[f: \langle R, \mathbb{T}_{GM} \rangle \rightarrow \langle \mathbb{R}^3, \mathbb{T}_{\mathbb{R}^3} \rangle\]
		$\mathbb{T}_{GM}$ provides a system for 3-dimensional region-based topology, in which regions and $\mathbb{R}^3$ are topologically equivalent. It results that $f$ is bijective, continuous, and has a continuous inverse.
	\end{enumerate}
	\section{Conclusion}
	In conclusion, motivated by the limitations identified in traditional qualitative spatial representation approaches particularly Goodman-style predicative mereologies and pseudo-topological frameworks, we have developed a unified theory bridging mereotopology and geometry by grounding topological reasoning in regular open sets and Tarski’s geometry of solids. By reconstructing points as higher-level abstractions over sets of regions and interpreting mereological classes as topological individuals, our approach constructively derives Tarski’s axioms within a purely nominal and inductive type-theoretic framework. This point-free formalism offers enhanced expressiveness, a securely founded and general theory with semantic precision, outperforming classical region-based calculi and first-order logic approaches. 
	\par Formalized within the Coq proof assistant using dependent type theory and classical logic, our system ensures rigorous spatial reasoning capabilities, certified proofs, and ontological robustness. This makes it particularly suitable for high-assurance applications where rigorous formalization of spatial propositions is crucial, including geographic information systems (GIS), spatial ontology modeling, architectural validation, and autonomous navigation tasks such as complex route planning and collision avoidance.
	\par Moreover, this theory should propose significant promise for symbolic supervision or fine-tuning of large language models (LLMs). The $\lambda$-MM library, enriched with logically structured datasets and spatial reasoning tasks, provides higher-quality training resources compared to traditional datasets. This would facilitate effective neuro-symbolic integration by offering symbolic scaffolds that enhance the spatial reasoning abilities of generative AI systems, bridging formal logical reasoning and generative modeling.
	\par Looking forward, our approach supports modular extensions and invites further integration of advanced spatial properties such as convexity and connectivity. Future research will explore these directions, aiming to further enhance representational power, computational applicability, and integration potential within hybrid neuro-symbolic learning pipelines.
	\bibliography{paper_20381_bib}
\end{document}